\documentclass{elsart}
\pdfoutput=1

\usepackage{graphicx}
\usepackage{amssymb, amsmath}
\usepackage{subfigure}

\begin{document}

\begin{frontmatter}



\title{Boundary and shape of binary images}


\author{Birgit van Dalen}

\address{Mathematisch Instituut, Universiteit Leiden, Niels Bohrweg 1, 2333 CA Leiden, The Netherlands}

\ead{dalen@math.leidenuniv.nl}

\begin{abstract}
In this paper we will consider an unknown binary image, of which the length of the boundary and the area of the image are given. These two values together contain some information about the general shape of the image. We will study two properties of the shape in particular. Firstly, we will prove sharp lower bounds for the size of the largest connected component. Secondly, we will derive some results about the size of the largest ball containing only ones, both in the case that the connected components are all simply connected and in the general case.
\end{abstract}

\end{frontmatter}

% main text
\section{Introduction}

Digital pictures, images consisting of pixels with discrete values, have been studied for several decades. The field of discrete geometry is concerned with the geometric properties of digital pictures, such as area, shape, boundary and connectivity \cite{kletterosenfeld,digtop}. This is used in, among other things, the segmentation, thinning and boundary-detection of images, which have applications in industry and medical imaging \cite{thinning,boundarydetection,segmentation}.

Digital pictures are usually 2- or 3-dimensional. Various grids can be used, such as square grids and hexagonal grids \cite{gray}. In this paper we consider binary digital pictures on a 2-dimensional square grid. The picture is a rectangle consisting of pixels or cells, i.e. unit squares that has value 0 or 1. The number of cells with value 1 is called the \emph{area} of the picture. Two cells are called \emph{4-adjacent} if they have an edge in common, and \emph{8-adjacent} if they have at least a vertex in common. The \emph{boundary} of a digital picture can be defined as the pairs consisting of two adjacent cells, one with value 0 and one with value 1 \cite{connectivity}. If we do this for 4-adjacency, then the boundary corresponds to the edges that separate the cells with value 1 from the cells with value 0. The number of such edges is called the \emph{length of the boundary} or sometimes the \emph{perimeter length} \cite{gray}.

In this paper we will consider an unknown binary image, of which the length of the boundary and the area of the picture are given. These two values together contain some information about the general shape of the picture. We will study two properties of the shape in particular. Firstly, using 4-adjacency, we can define the connected components of the picture \cite{connectivity}. We will prove sharp lower bounds for the size of the largest connected component.

The second question that we are interested in is: what is the size of the largest ball containing only ones? Or equivalently, considering for each cell the city block distance to the boundary \cite{distance}, what is the maximal distance that occurs? We will derive some results about this question, both in the case that the connected components are all simply connected (that is, they do not have any holes \cite{connectivity}) and in the general case.

After introducing some notation in Section \ref{notation}, we will tackle the first question in Section \ref{component} and the second question in Section \ref{ball}.

\section{Definitions and notation}\label{notation}

Let a \emph{cell} in $\mathbb{R}^2$ be a square of side length 1 of which the vertices have integer coordinates. A \emph{binary image} is a rectangle in $\mathbb{R}^2$ consisting of a number of cells, such that each cell inside the rectangle has been assigned a value 0 or 1. We will often refer to \emph{a one} or \emph{a zero} of a binary image, meaning a cell that has been assigned that value. When exactly $N$ of the cells of a binary image have been assigned the value 1, we say that the image \emph{consists of $N$ ones}.

We will only consider 4-adjacency \cite{connectivity}, hence we will simply call two cells \emph{neighbours} if they have a common edge. Two cells $c$ and $c'$ with value 1 in a binary image are called \emph{connected} if there is a path $c=c_1, c_2, \ldots, c_n = c'$ of cells with value 1 such that $c_i$ and $c_{i+1}$ are neighbours for $1 \leq i \leq n-1$. Being connected is an equivalence relation and the equivalence classes are called the \emph{connected components} of the image.

A connected component is said to contain a \emph{hole} if there is a zero or a group of zeroes that is completely surrounded by ones of the connected component.

\begin{figure}
  \begin{center}
    \subfigure[The length of the boundary of this image is 34.]{\label{figurebinary1}\includegraphics{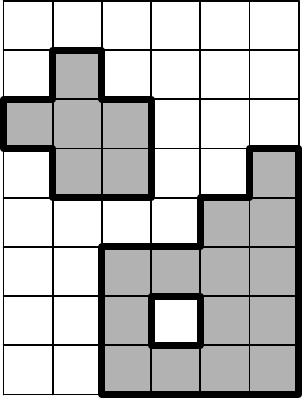}}
    \qquad
    \subfigure[In each cell with value 1 the distance to the boundary is indicated.]{\label{figurebinary2}\includegraphics{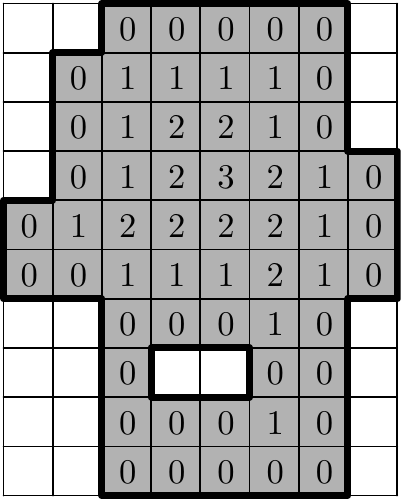}}
    \qquad
    \subfigure[A ball with radius 3.]{\label{figurebinary3}\includegraphics{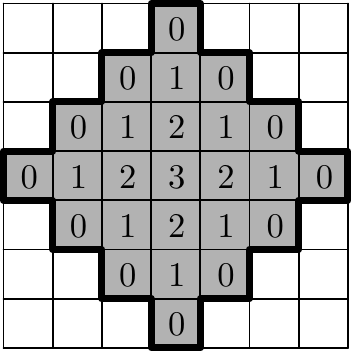}}
  \end{center}
\caption{Three binary images. The grey cells have value 1.}
% \label{blabla}
\end{figure}

The \emph{boundary} of a binary image consists of edges of cells. An edge belongs to the boundary if
\vspace{-\baselineskip}
\begin{itemize}
\item it is the common edge of two neighbouring cells, one of which has value 1 and one of which has value 0, or
\item it belongs to exactly one cell within the rectangle (i.e. it is part of the outer edge of the rectangle) and that cell has value 1.
\end{itemize}
We define the \emph{length of the boundary} as the number of edges that belong to the boundary. A binary image with its boundary is shown in Figure \ref{figurebinary1}.

For each cell $c$ with value 1 in a binary image, we define the \emph{distance to the boundary} $d(c)$ recursively. A cell of which one of the edges belongs to the boundary, has distance 0 to the boundary. For any other cell $c$ with value 1, we set
\[
d(c) = 1 + \min\{\ d(c') \ | \ c' \text{ and } c \text{ are neighbours }\}.
\]
See Figure \ref{figurebinary2} for an example. In the literature this specific distance function is often referred to as city block distance \cite{distance}.

For any integer $i \geq 1$ we define the \emph{$i$-boundary} similarly to the boundary. An edge belongs to the $i$-boundary if it is a common edge of two cells with value 1, one of which has distance $i-1$ to the boundary and the other of which has distance $i$ to the boundary. The $i$-boundary separates the cells $c$ with value 1 and $d(c) \geq i$ from the cells $c$ with value 0 or $d(c) \leq i-1$.

We say that a binary image contains a \emph{ball with radius $k$} if there is a cell with value 1 that has distance at least $k$ to the boundary. In that case the connected component containing this cell must contain at least $2k^2 + 2k + 1$ cells. See also Figure \ref{figurebinary3}.

\section{Largest connected component}\label{component}

Let $F$ be a binary image consisting of $m^2$ ones. If the ones are arranged into one square with side length $m$, then the boundary of $F$ has length $4m$. This is the smallest possible boundary for this number of ones (see also Lemma \ref{lemma4wortel}). If the length of the boundary is greater than $4m$, then the image may contain more than one connected component. We can, however, still prove a good lower bound on the size of the largest connected component. We will do this in two cases: when the boundary has length $4m$ plus some constant, and when the boundary has length $4m$ times some constant. In the second case we will also generalise to an image consisting of $N$ ones, where $N$ does not need to be a square.

First we prove two lemmas.

\begin{lem}\label{lemmaminimalvalue}
Let $r \geq 2$ and $0 \leq A < B$ be integers and let $S$ be an integer satisfying $rA \leq S \leq rB$.
The minimal value of
\[
f(k_1, k_2, \ldots, k_r) = \sqrt{k_1} + \sqrt{k_2} + \ldots + \sqrt{k_r}
\]
where $k_1, k_2, \ldots, k_r$ are integers in the interval $[A, B]$ for which $k_1 + k_2 + \cdots + k_r = S$, is attained at some $r$-tuple $(k_1, k_2, \ldots, k_r)$ for which $k_i \not\in \{A, B\}$ holds for at most one value of $i$.
\end{lem}

\begin{pf}
We argue by contradiction. Suppose the minimal value of $f$ is attained at some $r$-tuple $(k_1, k_2, \ldots, k_r)$ for which we have $k_1, k_2 \not\in \{A, B\}$. Let $S' = k_1+k_2$. Consider all possible values of $g(x) = \sqrt{x} + \sqrt{S'-x}$, where $x$ is an integer in the interval $[A, B]$ such that $S'-x \in [A,B]$ as well. Our assumption implies that the minimal value of $g$ is attained when $x=k_1$ and also when $x=k_2$. We now distinguish between two cases.

First suppose $k_1 + k_2 \leq A+B$. When we take $x=A$, we have $S'-x = k_1 + k_2 - A \leq B$ and $S'-x \geq A$, so $S'-x \in [A,B]$. Hence by our assumption $g(A) \geq g(k_1) = g(k_2)$. On the other hand, the continuous function $g(x) = \sqrt{x} + \sqrt{S'-x}$ on the interval $[0, S'] \subset \mathbb{R}$ is monotonically increasing on $[0, S'/2]$ and monotonically decreasing on $[S'/2, S']$. At least one of $k_1$, $k_2$ must be in $[0, S'/2]$ and $A < k_1, k_2$, so we must have $g(A) < g(k_1) = g(k_2)$, which yields a contradiction.

Now suppose $k_1 + k_2 > A+B$. When we take $x=B$, we have $S'-x = k_1 + k_2 - B > A$ and $S'-x \leq B$, so $S'-x \in [A,B]$. Similarly to above, this leads to a contradiction.
\end{pf}

\begin{lem}\label{lemma4wortel}
Let $k$ be a positive integer. A binary image consisting of $k$ ones has a boundary of length at least $4 \sqrt{k}$.
\end{lem}

\begin{pf}
First suppose that there is just one connected component. Let the smallest rectangle containing the component have side lengths $a$ and $b$. The boundary of the rectangle has length equal to or smaller than the boundary of the original image, so the boundary of the image has length at least $2a + 2b$. On the other hand, we have $k \leq ab$, since all $k$ ones are contained in the rectangle. As $\frac{a+b}{2} \geq \sqrt{ab} \geq \sqrt{k}$, the boundary has length at least $4 \sqrt{k}$.

Now suppose that there are $r$ connected components consisting of $k_1$, $k_2$, \ldots, $k_r$ ones respectively. Then the boundary of the image has length at least $4\sqrt{k_1} + 4\sqrt{k_2} + \cdots + 4\sqrt{k_r}$. So it suffices to prove
\[
\sqrt{k_1} + \sqrt{k_2} + \cdots + \sqrt{k_r} \geq \sqrt{k_1 + k_2 + \cdots + k_r},
\]
which can easily be done by squaring both sides.
\end{pf}

We will now prove our first theorem, concerning an image with boundary only an additive constant larger than the minimal length.

\begin{thm}
Let $m$ and $c$ be positive integers. Suppose a binary image $F$ consists of $m^2$ ones and has a boundary of length $4m+4c$. If $m$ is sufficiently large compared to $c$, then the largest connected component of $F$ consists of at least $m^2 - c^2$ ones.
\end{thm}

\begin{pf}
Suppose to the contrary that the largest connected component of $F$ consists of $t \leq m^2 - c^2 - 1$ ones. We distinguish between two cases. First assume that $t \geq c^2 + 1$. By Lemma \ref{lemma4wortel} the boundary has length at least $4 \sqrt{t} + 4 \sqrt{m^2 - t}$, while it is given to be equal to $4m+4c$. So we have
\[
\sqrt{t} + \sqrt{m^2 - t} \leq m+c.
\]
By Lemma \ref{lemmaminimalvalue} the smallest possible value of $\sqrt{t} + \sqrt{m^2-t}$ is attained when $t = m^2 - c^2 - 1$ (and when $t = c^2 + 1$). So we must have
\[
\sqrt{m^2 - c^2 -1} + \sqrt{c^2 + 1} \leq m+c.
\]
Subtracting $\sqrt{c^2+1}$ from both sides and squaring gives
\[
m^2 - c^2 - 1 \leq m^2 + 2mc + 2c^2 + 1 - 2(m+c)\sqrt{c^2+1}.
\]
This is equivalent to
\[
m \leq \frac{3c^2 + 2 - 2c\sqrt{c^2+1}}{2\sqrt{c^2+1} - 2c}.
\]
Hence for sufficiently large $m$, this case is impossible.

Now consider the case that $t \leq c^2$. Suppose we have $r$ connected components. Then $r \geq \frac{m^2}{t} \geq \frac{m^2}{c^2}$. The boundary of each connected component has length at least 4, so the total length of the boundary is at least $4r \geq 4 \frac{m^2}{c^2}$. Therefore, we must have
\[
\frac{m^2}{c^2} \leq m + c.
\]
For sufficiently large $m$, this is also impossible. We conclude that the largest connected component must consist of at least $m^2 - c^2$ ones.
\end{pf}

The bound given in this theorem is sharp: suppose the ones in the image are grouped in two connected components, an $(m-c) \times (m+c)$ rectangle and a $c \times c$ square. The boundary of the rectangle then has length $4m$, while the boundary of the square has length $4c$, so in total the boundary of $F$ has length $4m+4c$.

The next theorem concerns a binary image consisting of $m^2$ ones and having a boundary of length a constant times $4m$.

\begin{thm}\label{theorem4mc}
Let $m$ and $c$ be positive integers such that $m$ is divisible by $c$ and $m \geq c(c+1)$. Suppose a binary image $F$ consists of $m^2$ ones and has a boundary of length $4mc$. Then the largest connected component of $F$ consists of at least $\frac{m^2}{c^2}$ ones.
\end{thm}

\begin{pf}
Let $n$ be an integer such that $m = nc$. Then $F$ contains $c^2n^2$ ones and the boundary of $F$ has length $4c^2n$. We want to prove that the largest connected component of $F$ consists of at least $n^2$ ones. Suppose to the contrary that the largest connected component of $F$ consists of $t \leq n^2 - 1$ ones. Let $r$ be the number of connected components, and let $k_i$ be the number of ones in the $i$-th component, $1 \leq i \leq r$. Then by Lemma \ref{lemma4wortel} the boundary of $F$ is at least equal to
\begin{equation}\label{eqboundary}
4 \left( \sqrt{k_1} + \sqrt{k_2} + \cdots + \sqrt{k_r} \right).
\end{equation}
We will try to determine the minimal value of this and show that it is greater than $4c^2 n$.

The integers $k_1$, \ldots, $k_r$ are all in the interval $[1, t]$ and at least one of them is equal to $t$. For our purposes we may as well assume that $k_i \in [1, n^2 -1]$: by doing so we may find a minimal value that is even smaller than the actual minimal value, but if we can still prove that it is greater than $4c^2n$, we are done anyway.

The integers $k_1$, \ldots, $k_r$ furthermore satisfy $k_1 + k_2 + \cdots + k_r = c^2 n^2$. Also, since $c^2 \cdot (n^2-1) < c^2 n^2$, we know that $r \geq c^2 +1$.

By Lemma \ref{lemmaminimalvalue} the minimal value is attained at some $r$-tuple $(k_1, \ldots, k_r)$ of which at least $r-1$ elements are equal to 1 or $n^2-1$. Up to order, there is only one such $r$-tuple satisfying $k_1 + \cdots + k_r = c^2 n^2$. After all, suppose there are two such $r$-tuples, $(k_1 \leq k_2 \leq \ldots \leq k_r)$ and $(k_1' \leq k_2' \leq \ldots \leq k_r')$. Let $i$ be such that $k_i = 1$, $k_{i+1} > 1$ and let $j$ be such that $k_j' = 1$, $k_{j+1}' > 1$. If $i=j$, then the two $r$-tuples must be equal, as the sum of the elements is equal. So assume that $i \neq j$, say, $i > j$. Then $k_{i+2} = \ldots = k_r = n^2 -1$ and $k_{j+2}' = \ldots = k_r' = n^2 - 1$. Since the two sums of the $r$-tuples must be equal, we must have $k_{i+1} - k_{j+1}' = (i-j)(n^2-2)$. Since $k_{j+1}' \geq 2$ and $k_{i+1} \leq n^2-1$, the left-hand side can be at most $n^2 - 3$, while the right-hand side is at least $n^2 - 2$, which is a contradiction.

The unique $r$-tuple (ordered non-decreasingly) that satisfies the requirements is given by
\[
k_1 = \ldots = k_{r-v-1} = 1, \quad k_{r-v} = (c^2-v)n^2 + 2v + 1 - r, \quad k_{r-v+1} = \ldots = k_r = n^2 - 1,
\]
where $v$ is the unique positive integer such that
\[
(c^2-v-1)n^2 + 2v + 3 \leq r \leq (c^2-v)n^2 + 2v.
\]
This $r$-tuple must give the minimal value of (\ref{eqboundary}) under the conditions that $k_i \in [1, n^2-1]$ and $k_1 + \cdots + k_r = c^2n^2$. Therefore it now suffices to prove that
\begin{equation}\label{eqinequality}
(r-v-1) + \sqrt{(c^2-v)n^2 + 2v + 1 - r} + v \sqrt{n^2-1} > c^2n.
\end{equation}
From $m \geq c(c+1)$ we have $n \geq c+1$. This implies $n^2 > c^2 +1$, and from that we derive $v \leq c^2$: if $v \geq c^2+1$, then $\sum_i k_i \geq (c^2+1)(n^2-1)= c^2n^2+n^2-c^2-1 > c^2n^2$, which contradicts $\sum_i k_i = c^2n^2$. We now distinguish between two cases: $v \leq c^2 -1$ and $v = c^2$.

First suppose $v \leq c^2 - 1$. Consider the function $f(x) = x + \sqrt{S-x}$ on the interval $[A, S-1]$. Its derivative is $f'(x) = 1 - \frac{1}{2\sqrt{S-x}}$, which is positive for $x \leq S-1$, so the function is strictly increasing on the interval. Hence for all $x \in [A, S-1]$ we have $f(x) \geq f(A)$. If we apply this for $A =(c^2-v-1)n^2 + 2v + 3$, $S=(c^2-v)n^2 + 2v + 1$ and $x = r$, we find that
\[
(r-v-1) + \sqrt{(c^2-v)n^2 + 2v + 1 - r} \geq
(c^2-v-1)n^2 + v + 2 + \sqrt{n^2-2}.
\]
As $n \geq c+1 \geq 2$, we have $n^2 - 2 \geq (n-1)^2$, hence the left-hand side of (\ref{eqinequality}) is at least
\[
(c^2-v-1)n^2 + v + 2 + \sqrt{(n-1)^2} + v \sqrt{(n-1)^2}
\]
As $c^2 - v - 1 \geq 0$ and $n^2 \geq n$, this is at least
\[
(c^2-v-1)n + v + 2 + (v+1)(n-1) = c^2n + 1 > c^2n,
\]
which proves that (\ref{eqinequality}) holds in this case.

Now suppose $v = c^2$. Then $r \leq 2c^2$. Recall that we also have $r \geq c^2 +1$. We have to prove
\[
r-c^2-1 + \sqrt{2c^2 + 1 - r} + c^2 \sqrt{n^2-1} > c^2n.
\]
We again apply $f(x) \geq f(A)$ with $f(x)$ as above, now with $A = c^2 + 1$, $S = 2c^2+1$ and $x=r$. We find
\[
r-c^2 -1 + \sqrt{2c^2 + 1-r} \geq (c^2+1) -c^2 - 1 + \sqrt{2c^2 + 1 -(c^2+1)} = c.
\]
Hence it suffices to prove
\[
c + c^2 \sqrt{n^2-1} > c^2 n.
\]
This is equivalent to
\[
c^4(n^2-1) > (c^2 n - c)^2,
\]
which we can rewrite as
\[
n > \tfrac{1}{2}(c + \tfrac{1}{c}).
\]
This follows from $n \geq c+1$, hence (\ref{eqinequality}) holds in this case as well. This completes the proof of the theorem.
\end{pf}

The bound given in this theorem is sharp: suppose the ones in the image are grouped in $c^2$ squares of side length $\frac{m}{c}$, containing $\frac{m^2}{c^2}$ ones each. Then the boundary of each square has length $4\frac{m}{c}$, so in total the boundary of $F$ has length $4mc$.

The condition that $m$, $c$ and $\frac{m}{c}$ be integers does not seem to be very essential in the above theorem or proof. In fact, in a similar way (though slightly more technical) we can prove a more general result in which this condition is omitted.

\begin{thm}\label{theorem4mcgen}
Let $N$ be a positive integer and $c > 1$ a real number. Suppose a binary image $F$ consist of $N$ ones and has a boundary of length at most $4c\sqrt{N}$. If $N$ is sufficiently large compared to $c$, then the largest connected component of $F$ consists of more than $\frac{N}{c^2}-1$ ones.
\end{thm}

\begin{pf}
Let $q = \frac{\sqrt{N}}{c} \in \mathbb{R}$. Then $F$ contains $c^2q^2$ ones and the boundary has length at most $4c^2q$. Let $1 \leq \varepsilon < 2$ be such that $q^2 - \varepsilon$ is an integer, and suppose there are $t \leq q^2 - \varepsilon$ ones in the largest connected component of $F$. We will derive a contradiction, from which the theorem then follows. Let $r$ be the number of connected components, and let $k_i$ be the number of ones in the $i$-th connected component, $1 \leq i \leq r$.

Similarly to the proof of Theorem \ref{theorem4mc} it suffices to prove that (for sufficiently large $q$ compared to $c$) the minimal value of
\[
\sqrt{k_1} + \sqrt{k_2} + \cdots + \sqrt{k_r},
\]
where $k_1$, \ldots, $k_r$ are integers in the interval $[1, q^2 - \varepsilon]$ satisfying $k_1 + k_2 + \cdots + k_r = c^2 q^2$, is greater than $c^2q$. Also similarly to the proof of Theorem \ref{theorem4mc}, that minimal value is attained when
\[
k_1 = \ldots = k_{r-v-1} = 1, \quad k_{r-v} = (c^2-v)q^2 + (\varepsilon +1)v + 1 - r, \quad k_{r-v+1} = \ldots = k_r = q^2 - \varepsilon,
\]
where $v$ is the unique positive integer such that
\[
(c^2-v-1)q^2 + (\varepsilon + 1)v + \varepsilon + 2 \leq r \leq (c^2-v)q^2 + (\varepsilon + 1)v.
\]
It suffices to prove that
\begin{equation}\label{eqinequalitygen}
(r-v-1) + \sqrt{(c^2-v)q^2 + (\varepsilon +1)v + 1 - r} + v \sqrt{q^2-\varepsilon} > c^2 q.
\end{equation}
Let $c^2 + \delta$ be the smallest integer strictly greater than $c^2$. Then we can choose $q$ large enough such that $\delta q^2 > 2(c^2+\delta)$, which is equivalent with $(c^2+\delta)(q^2-2) > c^2q^2$. As $\varepsilon < 2$, we then also have $(c^2+\delta)(q^2 - \varepsilon) > c^2q^2$. As $c^2q^2 \geq v(q^2-\varepsilon)$, we find $v \leq c^2 + \delta - 1 \leq c^2.$ We now distinguish between three cases: the case $v \leq c^2 - 1$, the case $c^2 -1 < v < c^2$ and the case $v = c^2$. (Note that depending on whether $c^2$ is an integer, only one of the two latter cases may occur.)

First suppose $v \leq c^2 - 1$. We have $r \geq (c^2-v-1)q^2 + (\varepsilon + 1)v + \varepsilon + 2$ and therefore (similarly to the proof of Theorem \ref{theorem4mc})
\[
(r-v-1) + \sqrt{(c^2-v)q^2 + (\varepsilon +1)v + 1 - r} \geq (c^2-v-1)q^2 + \varepsilon v + \varepsilon + 1 + \sqrt{q^2 - \varepsilon - 1}.
\]
Furthermore, assuming $q \geq 2$ we have $\sqrt{q^2-\varepsilon} > q-\varepsilon$ and $\sqrt{q^2-\varepsilon-1} \geq q-\varepsilon - 1$, hence the left-hand side of (\ref{eqinequalitygen}) is strictly greater than
\[
(c^2-v-1)q^2 + \varepsilon v + \varepsilon + 1 + (q - \varepsilon - 1) + v (q-\varepsilon) =
(c^2-v-1)q^2 + (v+1) q.
\]
As $c^2 - v -1 \geq 0$ and $q^2 \geq q$, is this at least
\[
(c^2 - v - 1)q + (v+1)q = c^2q,
\]
which proves (\ref{eqinequalitygen}) in this case.

Now suppose $c^2-1 < v < c^2$. The largest connected component of $F$ contains less than $q^2$ ones, and $F$ contains $c^2q^2$ ones, hence the number of connected components is greater than $c^2$. That implies
\[
(r-v-1) + \sqrt{(c^2-v)q^2 + (\varepsilon +1)v + 1 - r} \geq c^2-v-1 + \sqrt{(c^2-v)q^2 + (\varepsilon +1)v + 1 - c^2}.
\]
We have $(\varepsilon +1) v - c^2 + 1 > 0$, hence
\[
\sqrt{(c^2-v)q^2 + (\varepsilon +1)v + 1 - c^2} > \sqrt{(c^2-v)q^2} = q \sqrt{c^2-v}.
\]
Also, $c^2 - v - 1 > 0$ and (as above) $\sqrt{q^2 - \varepsilon} > q- \varepsilon$. Therefore it suffices to prove
\[
q \sqrt{c^2-v} + v (q-\varepsilon) \geq c^2 q,
\]
which is equivalent to
\[
(\sqrt{c^2-v} - (c^2 - v)) q \geq \varepsilon v.
\]
As $\varepsilon \leq 2$, it also suffices to prove
\[
(\sqrt{c^2-v} - (c^2 - v)) q \geq 2 v.
\]
Since $0 < c^2 - v < 1$, we have $(\sqrt{c^2-v} - (c^2 - v)) > 0$. Now note that for a given $c$, there is at most one possible value for $v$ satisfying $c^2- 1 < v < c^2$, as $v$ is an integer. This value does not depend on $q$. Therefore we can choose $q$ large enough such that it satisfies
\[
(\sqrt{c^2-v} - (c^2 - v)) q \geq 2 v.
\]
Hence (\ref{eqinequalitygen}) holds for sufficiently large $q$.

Finally suppose $v = c^2$. In this case (\ref{eqinequalitygen}) transforms into
\[
(r-c^2-1) + \sqrt{(\varepsilon +1)c^2 + 1 - r} + c^2 \sqrt{q^2-\varepsilon} > c^2 q.
\]
As above, we have $r \geq c^2$, hence
\[
(r-c^2-1) + \sqrt{(\varepsilon +1)c^2 + 1 - r} \geq (c^2-c^2-1) + \sqrt{(\varepsilon +1)c^2 + 1 - c^2} = -1 + \sqrt{\varepsilon c^2 + 1}.
\]
As $\varepsilon \geq 1$, we have $\sqrt{\varepsilon c^2 +1 } > c$. Also, $\varepsilon \leq 2$. Therefore it suffices to prove
\[
-1 + c + c^2 \sqrt{q^2-2} > c^2 q.
\]
After some rewriting, this is equivalent to
\[
q (2c^3-2c^2) \geq 2c^4 + c^2 - 2c + 1.
\]
Since $2c^3 - 2c^2 > 0$, this is true for sufficiently large $q$. Hence also in this case (\ref{eqinequalitygen}) holds for sufficiently large $q$. This completes the proof of the theorem.
\end{pf}

\section{Balls of ones in the image}\label{ball}

In the previous section we proved bounds on the size of the largest connected component of an image. However, we are also interested in the shapes of such components. It seems likely that if the boundary is small compared to the number of ones, then there needs to be a large ball-shaped cluster of ones somewhere in the image. In this section we will prove lower bounds on the radius of such a ball.

First we prove some lemmas about the length of the $i$-boundary of an image.

\begin{lem}\label{lemmadistance1}
In a binary image, the length of the $1$-boundary is at most three times the length of the boundary.
\end{lem}

\begin{pf}
We can split the boundary into a number of simple, closed paths. (If there is more than one way to do this, we just pick one.) Let $\mathcal{P}$ be one of those paths, and denote its length by $L_0$. Let $S$ be the set of cells that have value 1 and have an edge in common with $\mathcal{P}$. Either the cells in $S$ are all on the outside of the path, or they are all on the inside of the path. Let $L_1$ be the number of edges of cells in $S$ that are part of the 1-boundary. (These edges do not necessarily form a simple, closed path.) We will prove a bound on $L_1$ in terms of $L_0$.

Consider all the pairs of edges of $\mathcal{P}$ having a vertex in common. There are three possible configurations, as shown in Figure \ref{figureconnections}. We call a pair of edges that form a straight line segment a \emph{straight connection}. The other two types we call \emph{corners}. A corner is of type I if both edges belong to the same cell with value 0; it is of type II if both edges belong to the same cell with value 1.

\begin{figure}
\begin{center}
\includegraphics{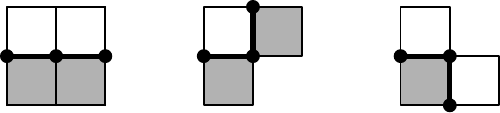}
\end{center}
\caption{From left to right: a straight connection, a corner of type I and a corner of type II. Such corners may also be called \emph{reentrant} and \emph{salient} respectively \cite{corners}.}
\label{figureconnections}
\end{figure}

We distinguish between three cases.

\textit{Case 1.} The path $\mathcal{P}$ consists of only four edges, and the cell enclosed by $\mathcal{P}$ has value 1. In this case $L_0 = 4$ and $L_1 = 0$.

\textit{Case 2.} The path $\mathcal{P}$ consists of more than four edges, and the cells in $S$ are on the inside of $\mathcal{P}$. Let $a$ be the number of straight connections and let $b$ be the number of corners of type I. Then the number of corners of type II must be $b+4$. We have $L_0 = a+2b+4$. Each edge of $\mathcal{P}$ is the edge of a cell in $S$, and each cell in $S$ has at least one edge in $\mathcal{P}$. In a corner of type II, we count the same cell in $S$ twice, so the number of cells in $S$ is $a+2b+4 - (b+4) = a+b$. Now we calculate an upper bound for $L_1$. Each cell in $S$ has four edges, of which in total $a+2b+4$ belong to $\mathcal{P}$. Also, the two cells in $S$ next to a straight connection share an edge that does not belong to either the boundary or the 1-boundary. Hence
\[
L_1 \leq 4(a+b) - (a+2b+4) - 2a = a+2b-4 = L_0 - 8.
\]

\textit{Case 3.} The cells in $S$ are on the outside of $\mathcal{P}$. Let $a$ be the number of straight connections and let $b$ be the number of corners of type I. Then $b \geq 4$ and there are $b-4$ corners of type II. Similarly to above, we find $L_0 = a + 2b - 4$, the number of cells in $S$ is $a+b$ and
\[
L_1 \leq 4(a+b) - (a+2b-4) - 2a = a+2b+4 = L_0 + 8.
\]
Since $L_0 \geq 4$, we have $L_1 \leq 3L_0$. This inequality obviously also holds in Cases 1 and 2.

Let $l_0$ be the length of the boundary and let $l_1$ be the length of the $1$-boundary of this image. Then $l_0$ is the sum of the lengths $L_0$ of all the paths $\mathcal{P}$, while $l_1$ is at most the sum of the lengths $L_1$ (we have counted each edge of the 1-boundary at least once). We conclude $l_1 \leq 3l_0$.
\end{pf}

\begin{lem}\label{lemmadistancei}
Let $i \geq 1$ be an integer. In a binary image, the length of the $(i+1)$-boundary is at most $\frac{2i+3}{2i+1}$ times the length of the $i$-boundary.
\end{lem}

\begin{pf}
Recall that the $i$-boundary consists of the edges between cells with distance $i-1$ to the boundary and cells with distance $i$ to the boundary. Just like the boundary, we can split the $i$-boundary into a number of simple, closed paths. Let $\mathcal{P}$ be one of those paths, and denote its length by $L_i$. Let $S$ be the set of cells that have distance $i$ to the boundary and have an edge in common with $\mathcal{P}$. Either the cells in $S$ are all on the outside of the path, or they are all on the inside of the path. Let $L_{i+1}$ be the number of edges of cells in $S$ that are part of the $(i+1)$-boundary. (These edges do not necessarily form a simple, closed path.) Analogously to the proof of Lemma \ref{lemmadistance1} we can prove a bound on $L_{i+1}$ in terms of $L_i$:
\begin{itemize}
\item In Case 1, $L_i = 4$ and $L_{i+1}=0$.
\item In Case 2, $L_{i+1} \leq L_i - 8$.
\item In Case 3, $L_{i+1} \leq L_i + 8$.
\end{itemize}
In Case 3, where in Lemma 3 we had $L_0 \geq 8$, we now have $L_i \geq 8i+4$. We will prove this here. Somewhere within $\mathcal{P}$ there must be a cell $c$ with value 0. A horizontal line drawn through $c$ must cross $\mathcal{P}$ somewhere to the left of $c$ and somewhere to the right of $c$. Between those two edges of $\mathcal{P}$ there must be at least $2i+1$ cells: $c$ and two cells at distance $j$ for each $j$ with $0 \leq j \leq i-1$. Similarly, there are at least $2i+1$ cells stacked in the vertical direction between two pieces of $\mathcal{P}$. Hence $L_i \geq 4(2i+1)$.

Since we have $L_{i+1} \leq L_i + 8$, we may conclude in Case 3 that
\[
\frac{L_{i+1}}{L_i} \leq 1 + \frac{8}{L_i} \leq 1 + \frac{8}{8i+4} = \frac{2i+3}{2i+1},
\]
and hence $L_{i+1} \leq \frac{2i+3}{2i+1} \cdot L_i$. Obviously this inequality holds in Cases 1 and 2 as well.

Let $l_{i}$ be the length of the $i$-boundary and let $l_{i+1}$ be the length of the $(i+1)$-boundary of this image. As in the proof of Lemma \ref{lemmadistance1} we conclude $l_{i+1} \leq \frac{2i+3}{2i+1} l_i$.
\end{pf}

\begin{lem}\label{lemmanumber}
Let $i \geq 0$ be an integer. In a binary image, the number of cells at distance $i$ from the boundary is at most $2i+1$ times the length of the boundary.
\end{lem}

\begin{pf}
For $i \geq 0$, let $A_i$ be the number of cells at distance $i$ from the boundary. For $i \geq 1$, let $l_i$ be the length of the $i$-boundary. Let $l_0$ be the length of the boundary. Each cell at distance $i$ from the boundary, $i \geq 1$, has at least one neighbour at distance $i-1$ from the boundary, hence the number of cells at distance $i$ from the boundary is at most equal to the length of the $i$-boundary. Similarly, the number of cells at distance 0 from the boundary is at most $l_0$. Furthermore, for $i \geq 1$ we have by Lemmas \ref{lemmadistance1} and \ref{lemmadistancei} that
\[
l_i \leq \frac{2i+1}{2i-1} \cdot l_{i-1} \leq \frac{2i+1}{2i-1} \cdot \frac{2i-1}{2i-3} \cdot l_{i-2} \leq \ldots \leq \frac{2i+1}{2i-1} \cdot \frac{2i-1}{2i-3} \cdot \cdots \cdot \frac{3}{1} \cdot l_0 = (2i+1)l_0.
\]
For $i=0$ it trivially holds that $l_i \leq (2i+1) l_0$. Hence for $i \geq 0$ we have
\[
A_i \leq (2i+1)l_0.
\]
\end{pf}

We now use these lemmas to prove our next theorem.

\begin{thm}\label{theoremball}
Let $N$ and $l$ be positive integers. Suppose a binary image $F$ consists of $N$ ones and has a boundary of length $l$. Then the image contains a ball of radius $\left\lceil \sqrt{\frac{N}{l}} -1 \right\rceil$.
\end{thm}

\begin{pf}
For $i \geq 0$, let $A_i$ be the number of cells with value 1 at distance $i$ from the boundary. Let $k$ be a positive integer. Recall that $F$ contains a ball with radius $k$ if there is a cell with value 1 that has distance at least $k$ to the boundary. Using Lemma \ref{lemmadistancei} we can find an upper bound for the number of cells with value 1 and distance to the boundary at most $k-1$:
\[
A_0 + A_1 + A_2 + \cdots + A_{k-1} \leq (1 + 3 + \cdots + 2k-1) l = k^2 l.
\]
Hence if $N > k^2 l$, then $F$ contains a ball with radius $k$.

Now let $k = \left\lceil \sqrt{\frac{N}{l}} -1 \right\rceil$ and assume that it is a positive integer (if it is not, then the theorem is trivial). Then $k < \sqrt{\frac{N}{l}}$, hence $N > k^2 l$. Therefore $F$ contains a ball with radius $\left\lceil \sqrt{\frac{N}{l}} -1 \right\rceil$.
\end{pf}

\begin{rem}
Suppose as in Theorem \ref{theorem4mcgen} that the boundary of $F$ has length $4c \sqrt{N}$ for some $c \in \mathbb{R}$. Then Theorem \ref{theoremball} says that $F$ contains a ball of radius $\left\lceil \sqrt{\frac{\sqrt{N}}{4c}} -1 \right\rceil$. This ball contains approximately $\frac{\sqrt{N}}{2c}$ ones. On the other hand, Theorem \ref{theorem4mcgen} tells us that there exists a connected component with more than $\frac{N}{c^2}-1$ ones. This is roughly four times the square of the size of the ball, but this component does not need to be ball-shaped.
\end{rem}

If the binary image contains no holes, then we can prove a much stronger result, by sharpening the lemmas in this section.

\begin{thm}\label{theoremnoholes}
Let $N$ and $l$ be positive integers. Suppose a binary image $F$ consists of $N$ ones and has a boundary of length $l$. Furthermore assume that none of the connected components of $F$ contains any holes. Then the image contains a ball of radius $\left\lfloor \frac{N}{l} \right\rfloor$.
\end{thm}

\begin{pf}
For $i \geq 0$, let $A_i$ be the number of cells with value 1 at distance $i$ from the boundary. Case 3 in the proofs of Lemmas \ref{lemmadistance1} and \ref{lemmadistancei} does not occur if the connected components of $F$ do not contain any holes. This means that in Lemma \ref{lemmadistance1} we can conclude that the length of the 1-boundary is strictly smaller than the length of the boundary, and in Lemma \ref{lemmadistancei} that the length of the $(i+1)$-boundary is strictly smaller than the length of the $i$-boundary. Hence we have for all $i \geq 0$
\[
A_i < A_{i-1} < \ldots < A_0 < l.
\]
Let $k$ be a positive integer. Then the number of cells with value 1 and distance to the boundary at most $k-1$ is
\[
A_0 + A_1 + A_2 + \cdots + A_{k-1} < kl.
\]
Hence if $N \geq kl$, then $F$ contains a ball of radius $k$. This is obviously the case for $k = \left\lfloor \frac{N}{l} \right\rfloor$.
\end{pf}

We will show by two examples that the bounds from the previous two theorems are nearly sharp.

\begin{figure}
  \begin{center}
    \subfigure[The binary image $F$ from the example, where $u=3$ and $c=2$.]{\label{figureholes1}\includegraphics[width=4cm]{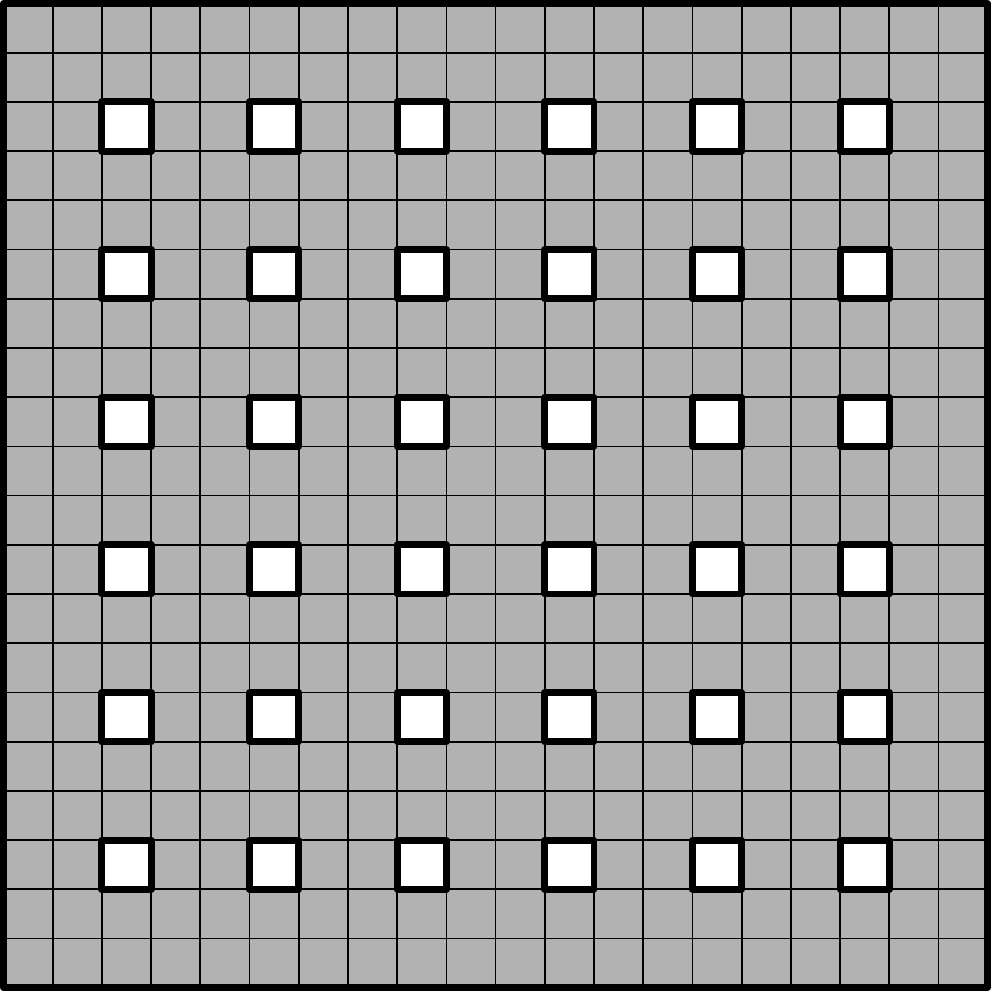}}
    \qquad
    \subfigure[When $u$ is odd, the radius of the largest ball that fits in the image is $u-2$.]{\label{figureholes2}\includegraphics{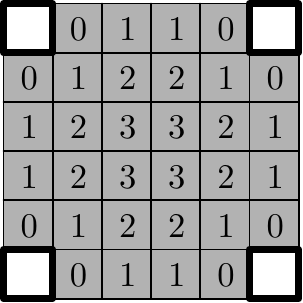}}
    \qquad
    \subfigure[When $u$ is even, the radius of the largest ball that fits in the image is $u-1$.]{\label{figureholes3}\includegraphics{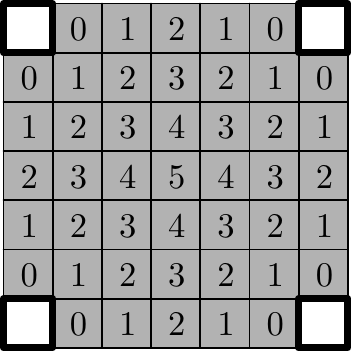}}
  \end{center}
  \caption{Some illustrations for Example \ref{exampleholes}.}
  \label{figureholes}
\end{figure}

\begin{exmp}\label{exampleholes}
Let $u$ and $c$ be positive integers. Consider a square of ones of side length $cu^2 + u -1$. Denote the cells in the square by coordinates $(i,j)$, where $1 \leq i,j \leq cu^2 + u - 1$. For all $i$ and $j$ that are divisible by $u$, we change the value of cell $(i,j)$ from 1 to 0. Let $F$ be the resulting binary image (see also Figure \ref{figureholes1}). The number of ones of $F$ is
\[
N = (cu^2 + u - 1)^2 - (cu)^2 = c^2u^4 + 2cu^3 + (- c^2-2c+1)u^2 - 2u + 1.
\]
The length of the boundary is
\[
l = 4(cu^2 + u -1 ) + 4c^2u^2 = 4(c^2+c)u^2 + 4u - 4.
\]
If $u$ is very large, we have $N \approx c^2u^4$ and $l \approx 4(c^2+2)u^2$. So according to Theorem \ref{theoremball}, $F$ should contain a ball of radius approximately
\[
\sqrt{\frac{N}{l}} \sim \sqrt{ \frac{c^2u^4}{4(c^2+c)u^2}} = \frac{1}{2} \cdot \sqrt{\frac{c^2}{c^2+c}} \cdot u, \qquad u \rightarrow \infty.
\]
If $u$ is odd, $F$ in fact contains a ball of radius $u-2$. If $u$ is even, then $F$ contains a ball of radius $u-1$. See also Figures \ref{figureholes2} and \ref{figureholes3}.
\end{exmp}

\begin{exmp}
Let $F$ consist of a rectangle of ones, with side lengths $a$ and $ta$, where $t \geq 1$. Then the number of ones is equal to $ta^2$, while the length of the boundary is equal to $2(t+1)a$. So according to Theorem \ref{theoremnoholes}, $F$ should contain a ball of radius $\lfloor\frac{ta^2}{2(t+1)a}\rfloor = \lfloor\frac{t}{t+1}\frac{a}{2}\rfloor$. The actual radius of the largest ball contained in $F$ is equal to $\left\lfloor \frac{a-1}{2} \right\rfloor$.
\end{exmp}

\end{document}